\documentclass[final]{siamltex}
\usepackage{graphics,graphicx,verbatim,textcomp,upquote}
\newlength\myverbindent 
\makeatletter
\def\verbatim@processline{%
 \hspace{\myverbindent}\the\verbatim@line\par}
\makeatother
\def\Im{\hbox{\rm\kern .3pt Im\kern 1pt}}
\def\pput(#1,#2)#3{\noindent\smash{\raise#2pt\hbox to 0pt
   {\kern #1pt #3\hss}}\ignorespaces}

\title{Exactness of quadrature formulas}

\author{Lloyd N.~Trefethen\thanks{\texttt{trefethen@maths.ox.ac.uk},
Mathematical Institute, University of Oxford, Oxford, OX2 6GG, UK.}}

\begin{document}

\maketitle


\begin{abstract}
The standard design principle for quadrature formulas is
that they should be exact for integrands of a given class,
such as polynomials of a fixed degree.  We show how
this principle fails to predict the actual 
behavior in four cases: Newton--Cotes, Clenshaw--Curtis, Gauss--Legendre, and
Gauss--Hermite quadrature.  Three further examples
are mentioned more briefly.
\end{abstract}


\pagestyle{myheadings}
\thispagestyle{plain}
\markboth{\sc Trefethen}
{\sc Exactness of quadrature formulas}

\section{\label{secintro}Introduction}
A quadrature formula is an approximation
\begin{equation}
I_n(f) = \sum_{j=1}^n w_j f(x_j)
\label{In}
\end{equation}
to a definite integral
\begin{equation}
I(f) = \int_D f(x)\kern .5pt dx.
\label{I}
\end{equation}
Here $D$ is a domain such as an interval in one dimension or a
hypercube in $s$ dimensions, the points $\{x_k\}$ are distinct
{\em nodes} in $D$, and the numbers $\{w_k\}$ are {\em weights.}
Sometimes a further weight function $w(x)$ is introduced in
(\ref{I}), as we shall see in section 5.\ \ Quadrature formulas
generally come in families defined by a rule that specifies how
the nodes and weights are determined for each choice of $n$, and we
shall use the word ``formula'' to refer to both the fixed $n$ case
and the family.  The aim with any quadrature formula is that the
error
\begin{equation}
E_n(f) = I_n(f)-I(f)
\label{error}
\end{equation}
should be small, and in particular,
one would like $E_n(f)$ to decrease rapidly as $n\to\infty$
when $f$ is smooth.

There is a standard design principle used in deriving quadrature
formulas: the formula should be exact when applied to a certain
class of integrands $f$.  In the case of quadrature on an interval,
this is often the set $P_{n-1}$ of polynomials of degree at most
$n-1$, in which case the result returned by the quadrature formula
is equal to the integral of the unique degree $n-1$ polynomial
interpolant through the data $\{f(x_k)\}$ at the points $\{x_k\}$.
This {\em exactness principle} has proved effective for a wide
range of problems.  Nevertheless, it is not a reliable guide to
the actual accuracy of quadrature formulas, as we shall show in
this article by considering four cases from this point of view:\
the Newton--Cotes, Clenshaw--Curtis, and Gauss formulas on
$[-1,1]$, and the Gauss--Hermite formula on $(-\infty,\infty)$.
The failure of the exactness principle is particularly extreme
in the cases of Newton--Cotes quadrature (as is well known)
and Gauss--Hermite quadrature (not so well known).

There are a number of excellent books on quadrature, including
the classic by Davis and Rabinowitz~\cite{dr} and the more recent
work by Brass and Petras~\cite{bp}.

\section{\label{secnc}Newton--Cotes quadrature}
Here and in the next two sections, our domain is the interval $D =
[-1,1]$.  The {\em Newton--Cotes formula,} going back to Isaac
Newton in 1676 and Roger Cotes in 1722, is the formula that
results from taking $\{x_k\}$ as equally spaced points from $-1$
to $1$, with $\{w_k\}$ determined so that $I_n(f) = I(f)$ for all
$f\in P_{n-1}$.  Most numerical analysis textbooks have a chapter
on numerical integration in which they discuss two quadrature
formulas.  First, the Newton--Cotes formula is introduced and
it is observed that it has polynomial exactness degree $n-1$.
Then Gauss quadrature is presented, based on optimal
points $\{x_k\}$ as defined by exactness degree (section~4),
and it is said to be better because it has exactness degree $2n-1$.

This is spectacularly misleading.  Gauss quadrature is indeed
better than Newton--Cotes, but this has nothing to do with its
doubled polynomial exactness degree.  In fact Clenshaw--Curtis
quadrature, as we shall discuss in the next section, has similar
behavior to Gauss quadrature but only the same exactness degree
$n-1$ as Newton--Cotes.  What makes the difference spectacular is
that Newton--Cotes doesn't merely converge less quickly as
$n\to\infty\kern .5pt$;
it diverges at an exponential rate, even for many analytic
integrands.  Meanwhile the Gauss and Clenshaw--Curtis formulas
converge for all continuous $f$, and at an exponential rate if $f$
is analytic.

\begin{figure}
\begin{center}
\vskip 15pt
\includegraphics[scale=.8]{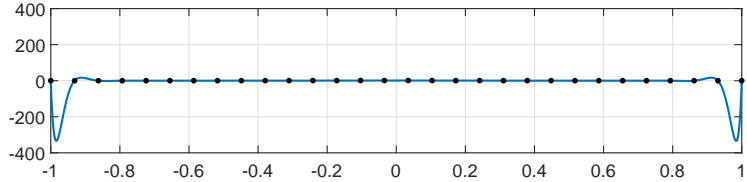}
\end{center}
\caption{\label{fignc1}Degree $29$ polynomial interpolant
to the integrand $f(x) = 1/(1+25x^2)$ in $30$ equispaced points
of $[-1,1]$.  The Runge phenomenon of oscillations near the
boundary leads to exponential divergence of the Newton--Cotes formula
as $n\to\infty$.}
\end{figure}

The failure of Newton--Cotes quadrature for larger values of
$n$ is well known.  This became apparent
to experts after the appearance in 1901 of Runge's paper on
polynomial interpolation in equispaced points~\cite{runge},
which shows that these interpolants experience oscillations whose
amplitude grows exponentially as $n\to\infty$, even for many
analytic functions $f$ (Figure~\ref{fignc1}).  The conclusion was
made rigorous by P\'olya in 1933~\cite{polya}.  P\'olya also showed
that convergence as $n\to\infty$ for all $f\in C([-1,1])$
occurs if and only if the sum $\sum_{j=1}^n |w_j|$ is bounded
as $n\to\infty$.  Since $\sum_{k=1}^n w_j = 2$, this condition
certainly holds if the weights are positive, as is the case
for the Gauss and Clenshaw--Curtis formulas.  The Newton--Cotes
weights, however, have alternating signs and grow in amplitude
at a rate of order $2^n$ as $n\to\infty$.\footnote{See the final
formula of~\cite{ousp}, from 1925, after which Ouspensky writes
``One sees that the coefficients $A_2, A_3, \dots, A_{n-2}$ tend
to infinity, making it evident that the Cotes formula loses all
practical value as the number of ordinates grows considerable.''}
Figure~\ref{fignc1} illustrates the divergence of the Newton--Cotes
formula with a plot of the degree $29$ polynomial interpolant
to $f(x) = 1/(1+25 x^2)$ (the celebrated Runge function).  The
corresponding quadrature estimate is $I_{30}(f) \approx -21.8$.
Even the sign is wrong, and the amplitude grows exponentially,
with $I_{50}(f) \approx -24{,}965$, for example.

\begin{table}
\caption{\label{tabnc1}Errors in Newton--Cotes integration of
polynomials $x^k$ and $T_k(x)$ with $n=30$. (Only even values of\/ $k$
are shown since the integrals and errors are zero when $k$ is odd.)
The errors for $T_k(x)$ are huge for $k\ge n$, a reflection of exponentially
large quadrature weights of alternating signs.}
\begin{center}
\begin{tabular}{c c c}
$~~~k~~~$ & $|E_n(x^k)|$ & $|E_n(T_k(x))|$ \\[2pt]
\hline
\vdots & \vdots & \vdots \\
26 & 0 & 0\\
28 & 0 & 0\\
\hline
30\vrule height 11pt width 0pt & .0000007 & 399.5 \\
32 & .000005 & 2711.1 \\
34 & .00002 & 8923.1 \\
36 & .00004 & 18765.9 \\
38 & .00009 & 27812.9 \\[2pt]
\smash{\vdots}\vphantom{0} & \smash{\vdots} & \smash{\vdots} 
\end{tabular}
\end{center}
\end{table}

What is not well known is how this failure relates to the
exactness principle, which we now consider in Table~\ref{tabnc1}.
With $n=30$, the Newton--Cotes formula integrates $x^k$ exactly for
$k< 30$, with $E_n(x^k)=0$.  A check of $E_n(x^k)$ for $k\ge 30$
looks unexpectedly promising, with the errors coming out very small,
smaller than one would have dreamed of for a quadrature formula
of exactness degree $n-1$.  However, this apparent good behavior
is an illusion associated with the exponential ill-conditioning
of monomial bases on $[-1,1]$.  (Numerically speaking, $x^k$ has
degree only $O(k^{1/2})$ for large $k$~\cite{nr}.)  Switching to
the well-conditioned basis of Chebyshev polynomials $T_k(x)$
reveals that huge errors set in immediately at degree $k=30$.

Evidently for the Newton--Cotes formula, exact integration
of degree $n-1$ polynomials has told us next to nothing about
accuracy in integrating other functions.

\section{\label{seccc}Clenshaw--Curtis quadrature}
The {\em Clenshaw--Curtis formula}, originating in
1960~\cite{cc}, consists of integrating the degree $n-1$ polynomial
interpolant through $n$ {\em Chebyshev points}
\begin{equation}
x_j = \cos(\kern .7pt j\pi/(n-1)), \quad 0\le j \le n-1.
\label{chebpts}
\end{equation}
This is the natural formula to apply in the context of
Chebyshev spectral collocation methods for differential equations,
and it is essentially the method by which Chebfun integrates a function, after
first reducing it to a polynomial of sufficiently high 
degree~\cite{chebfun}.  Alternatively, nodes and weights can
be computed explicitly in $O(n\log n)$ operations~\cite{wald}
and are available in Chebfun with the command {\tt [x,w] = chebpts(n)}.

Clenshaw--Curtis quadrature, like Newton--Cotes, has polynomial
exactness degree $n-1$, but with none of the misbehavior as
$n\to\infty$ since the weights are always positive.  P\'olya's
theory guarantees convergence for all $f\in C([-1,1])$ at a
rate that follows the smoothness of $f$; if $f$ is analytic, the
convergence is exponential~\cite[Theorems 19.3 and 19.4]{atap}.
Thus the obvious expectation for Clenshaw--Curtis quadrature is
that it should converge like Gauss quadrature, but at approximately
half the rate, since Gauss has polynomial exactness degree $2n-1$.

\begin{figure}
\begin{center}
\vskip 15pt
\includegraphics[scale=.8]{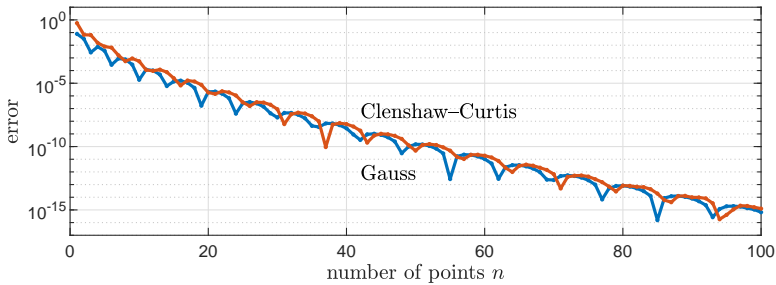}
\end{center}
\caption{\label{figcc1}Convergence of the Clenshaw--Curtis
and Gauss quadrature formulas for the integrand
$f(x) = \exp(-1/x^2)$.  The Gauss convergence is at
the expected rate; the surprise is that Clenshaw--Curtis
converges at this rate too.}
\end{figure}

This is not what happens.  Gauss quadrature behaves as expected;
the surprise is that Clenshaw--Curtis often converges at the
Gauss rate too~\cite[Theorem 19.5]{atap}.  Figure~\ref{figcc1}
shows a typical example.  This effect was noted experimentally
by Clenshaw and Curtis themselves, who wrote:\ ``We see that the
Chebyshev formula, which is much more convenient than the Gauss,
may sometimes nevertheless be of comparable accuracy''~\cite{cc}.
A paper on the subject was published by O'Hara and Smith in 1968,
who wrote: ``The Clenshaw--Curtis method gives results nearly
as accurate and Gauss\-ian quadratures for the same number of
abscissae''~\cite{ohs}.  Subsequently the effect was mentioned in
books of Evans~\cite{evans} and Kythe and Sch\"aferkotter~\cite{ks}
and then became more widely known through a paper of mine in
2008~\cite{gausscc}.

\begin{figure}
\vskip -10pt
\begin{center}
\includegraphics[scale=.8]{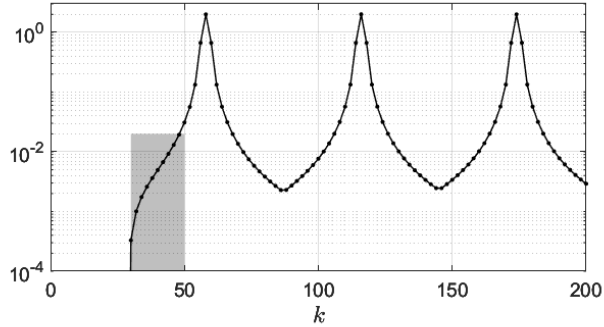}
\vspace{-110pt}
\end{center}
\caption{\label{figcc2}Errors in Clenshaw--Curtis integration
with $n=30$ of $T_k(x)$ as a function of even indices $k$.  The shading
highlights the small errors for $n\le k \ll 2n$.}
\end{figure}

O'Hara and Smith's explanation of the unexpected accuracy of
the Clenshaw--Curtis formula is that although its errors in
integrating degree $k$ polynomials are nonzero for $k\ge n$,
they are still very small for $n \le k \ll 2n$.  Note that this
is precisely a failure of exactness as a guide to accuracy.
The small errors are shown numerically in Table~\ref{tabcc1},
a repetition of Table~\ref{tabnc1} (without the distracting $x^k$
column) for Clenshaw--Curtis.  Figure~\ref{figcc2} gives a visual
picture of what is going on.

\begin{table}
\caption{\label{tabcc1}Errors in Clenshaw--Curtis integration of
polynomials $T_k(x)$ with $n=30$.  Though nonzero, the
errors for $n\le k \ll 2n$ are small, of order $O(n^{-2})$.  This
is the explanation given by O'Hara and Smith\/~{\rm \cite{ohs}} of the
unexpected accuracy of Clenshaw--Curtis.}
\begin{center}
\begin{tabular}{c c}
$~~~k~~~$ & $|E_n(T_k(x))|$ \\[2pt]
\hline
\vdots & \vdots \\
26 & 0 \\
28 & 0 \\
\hline
30\vrule height 11pt width 0pt & 0.0003 \\
32 & 0.001 \\
34 & 0.002 \\[2pt]
\smash{\vdots}\vphantom{0} & \smash{\vdots} \\
54 & 0.1 \\
56 & 0.7 \\
58 & 2.0 \\
\hline
60\vrule height 11pt width 0pt & 0.7 \\[2pt]
\smash{\vdots}\vphantom{0} & \smash{\vdots} 
\end{tabular}
\end{center}
\end{table}

The effect shown in Table~\ref{tabcc1} and
Figure~\ref{figcc2} leads readily to an understanding of 
the surprising convergence rate of Clenshaw--Curtis
quadrature as seen in Figure~\ref{figcc1}.
Any Lipschitz continuous integrand $f$ will have an absolutely
and uniformly convergent Chebyshev series,
\begin{equation}
f(x) = \sum_{j=0}^\infty a_j T_j(x),
\end{equation}
from which it follows that the error (\ref{error})  in Clenshaw--Curtis
quadrature is
\begin{equation}
E_n(f) = \sum_{\scriptstyle{j=n}\atop
\hbox{\scriptsize~$j$ even\kern 4pt}}^\infty a_j E_n(T_j).
\end{equation}
If all the errors $E_n(T_j)$ were of the same size $O(1)$,
then $E_n(f)$ would depend just on the Chebyshev coefficients
$\{a_j\}$, and this is approximately what happens in the regime
$2n\le j<\infty$.  For the ``shaded coefficients'' $a_j$ with
$n\le j \ll 2n$, however, the errors $E_n(T_j)$ are of size
$O(n^{-2})$, and whether they contribute significantly to $E_n(f)$
depends on the rate of decrease of $\{a_j\}$ as $j\to \infty$,
hence on the regularity of~$f$.  If $f$ is not analytic, as in
the example of Figure~\ref{figcc1}, then $\{a_j\}$ decrease slowly
enough that $a_jE_n(T_j)$ is much smaller for
$n \le j \ll 2n$ than for $j\approx 2n$,
making the contributions of the shaded coefficients
to $E_n(f)$ negligible and producing the doubled-degree effect in
its cleanest form.  (There is still a complication, however, in
that both Clenshaw--Curtis and Gauss converge at a rate faster than
expected by one power of $n$~\cite{xb}.)  If $f$ is analytic, then
$\{a_j\}$ decrease exponentially, and the terms $a_jE_n(T_j)$ for
$n \le j \ll 2n$ are no longer negligible
in comparison to the later ones, making the asymptotic
convergence rate of Clenshaw--Curtis indeed half that of Gauss.
For details, including the ``kink phenomenon'' observed for
Clenshaw--Curtis quadrature of
analytic integrands as the initial Gauss convergence rate cuts
in half after a certain value of $n$, see~\cite{gausscc,kink}.

\section{\label{secgauss}Gauss quadrature}
Gauss quadrature, discovered by Gauss in 1814~\cite{gauss},
is defined by having the maximum possible polynomial degree
of exactness, $2n-1$.  This is achieved by taking the nodes
$\{x_j\}$ as the roots of the degree $n$ Legendre polynomial
$P_n$, and the name {\em Gauss--Legendre} is also
used to distinguish this case from that of integrals in which
a nonconstant weight function is introduced in (\ref{I}) (next
section).  The weights for Gauss quadrature $\{w_j\}$ are all
positive, and the formula is extremely effective in practice.
Thanks to new algorithms introduced in the past 15 years and
available with the Chebfun command {\tt [x,w] = legpts(n)}, the
nodes and weights can be computed in a fraction of a second even
when $n$ is in the millions~\cite{bog,chebfun,haletown}.

Gauss quadrature is often described as optimal, but this is
only precisely true in senses that are tied to polynomials.
Specifically, it can be shown to be optimal by certain
measures for integrating functions that are analytic in a {\em
Bernstein ellipse,} an ellipse in the complex plane with foci
$\pm 1$; see~\cite{petras} and sections~4.9 and~6.9
of~\cite{bp}.  Analyticity in
an ellipse, however, is a skewed form of smoothness, requiring
more of a function in the middle of the interval than near
the endpoints.  Polynomials can resolve much faster wiggles
near endpoints than in the interior~\cite{dt}, and Gauss quadrature
(likewise Clenshaw--Curtis) inherits this property---as one sees
intuitively from its strong clustering of sample points near
the ends.  For example, Gauss quadrature converges faster for the
integrand $\sqrt{1.01-x}$ than for $\sqrt{0.1i - x}$, though the
singularity of the first function is ten times closer to $[-1,1]$
than that of the second.\footnote{In Chebfun, try {\tt cheb.x,
plotcoeffs([sqrt(1.01-x) sqrt(0.1i-x)])}.}

\begin{figure}
\begin{center}
\vskip 15pt
\includegraphics[scale=.8]{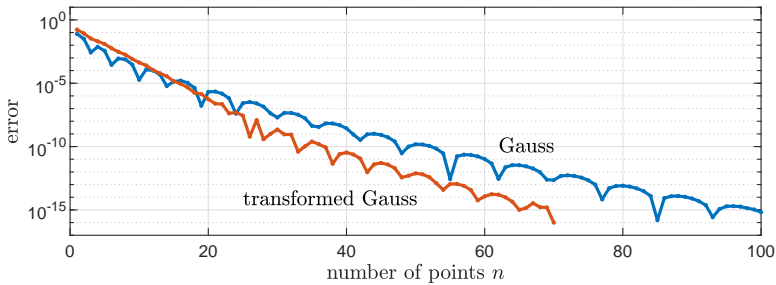}
\end{center}
\caption{\label{figtrans}Convergence of the 
Gauss formula and a conformally transformed Gauss formula $(\ref{Ing})$
for the same integrand
$f(x) = \exp(-1/x^2)$ as in Figure~$\ref{figcc1}$.  While
such transformed formulas have a limited role in practice,
they illustrate that Gauss quadrature is not
optimal.}
\end{figure}

From a user's point of view, it would seem more natural to
consider integrands with uniform smoothness across $[-1,1]$.
In the analytic case, one might require analyticity in an
$\varepsilon$-neighborhood of this interval, and quadrature
formulas based on this assumption can be derived by transplanting
Gauss quadrature via a conformal map $g$ with
$g([-1,1]) = [-1,1]$ of a Bernstein ellipse onto a cigar, or
more simply, an infinite strip.  Here the integral (\ref{I}) becomes
\begin{equation}
I(f) = \int_{-1}^1 g'(s) f(g(s))\kern .5pt ds,
\label{Ig}
\end{equation}
and applying Gauss quadrature in the $s$ variable gives
the transformed quadrature formula
\begin{equation}
I_n(f) = \sum_{j=1}^n w_j \kern .8pt g'(s_j) f(g(s_j)).
\label{Ing}
\end{equation}
An example is shown in Figure~\ref{figtrans}, with $g$ taken as
the conformal map of the Bernstein ellipse with parameter $\rho
= 1.4$ (the sum of the semiminor and semimajor axes) onto an
infinite strip.  Such transplanted formulas were introduced
in~\cite{hale}, and the conformal mapping idea goes back in the
theoretical literature to Bakhvalov in 1967~\cite{bakh} and was
applied for spectral methods by Kosloff and Tal-Ezer~\cite{kte}.
As expected, the transformed quadrature nodes are much more
uniformly distributed, with density ${\approx}\kern .7pt\pi/2$
times greater in the middle of the interval than for Gauss
quadrature (not shown).  In the limit $n\to\infty$ and
$\rho\downarrow 1$, $\approx$ becomes $\sim\kern .7pt$.

\begin{figure}
\begin{center}
\vskip 15pt
\includegraphics[scale=.8]{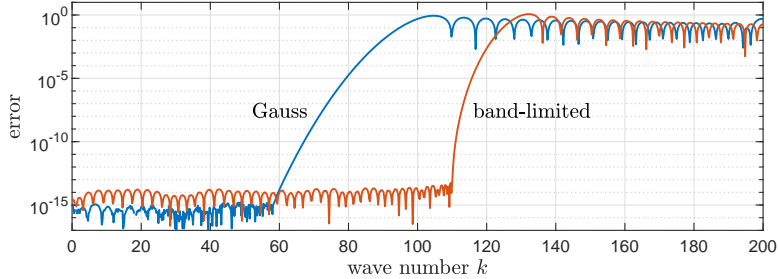}
\end{center}
\caption{\label{ggqcos} Errors of $50$-point Gauss and band-limited
quadrature formula integrating
$\exp(i k x)$ for $k\in [-c,c\kern .5pt ]$,
with $c$ given by $(\ref{c})$.
By this measure of integrating complex exponentials, the
band-limited formula is notably more efficient.}
\end{figure}

\begin{figure}
\begin{center}
\vskip 15pt
\includegraphics[scale=.8]{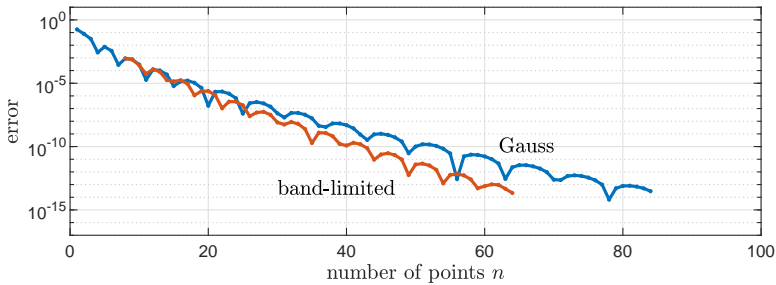}
\end{center}
\caption{\label{figtransggq}Comparison of
the Gauss and band-limited quadrature formulas
for the integrand $f(x) = \exp(-1/x^2)$, as in
Figures~$\ref{figcc1}$ and~$\ref{figtrans}$.  The
curves stop at accuracy $10^{-13}$ because our
computation of nodes and weights in the band-limited case
is a few digits off machine accuracy, as was
evident in Figure~$\ref{ggqcos}$.}
\end{figure}

A different approach to developing quadrature formulas with
more uniform behavior, {\em band-limited quadrature,}
is based on time-frequency analysis.  
Suppose one seeks a formula (\ref{In}) that will integrate the
functions $\exp(ikx)$ to high accuracy for all the wave numbers
$k\in [-c,c\kern .5pt]$ for some
$c>0$.  Note that this is a continuum of
wave numbers, not just the integer multiples of $\pi$ that would
make the integrand 2-periodic.  No choice of nodes and weights can
integrate the continuum exactly, because its dimension is infinite.
However, it has been known since the work of Slepian, Landau,
and Pollak at Bell Labs in the 1950\kern .3pt s and
1960\kern .3pt s that the 
numerical dimension of this space is finite, just a bit larger
than $2c/\pi$~\cite{prolateII,prolateIII,prolateI}.  Specifically,
the singular values $\sigma_j$ of the bivariate kernel
function $\exp(ikx)$ for $x\in[-1,1]$, $k\in[-c,c\kern .5pt]$
decrease exponentially for $j > 2\kern .5pt c/\pi$.\footnote{In
Chebfun, try {\tt K = chebfun2(@(x,k) exp(1i*k.*x),[-1 1 -c
c]),} {\tt semilogy(svd(K),\textquotesingle.\textquotesingle)}.}
By applying the method known as generalized Gauss
quadrature~\cite{bgr,orx} to an appropriate set of prolate spheroidal
wave functions (PSWFs), one can develop quadrature rules
that integrate these functions to high accuracy,
as shown by Xiao, Rokhlin, and Yarvin~\cite{orx,xry}.\footnote{Recent
experiments by Jim Bremer (unpublished) show that essentially the same results
can be obtained without use of PSWFs by applying generalized
Gauss quadrature directly to the functions $\exp(ikx)$.}
Later work has led to
fast methods for calculating the nodes and weights
of these formulas numerically~\cite{or,orx}.
For $15$-digit accuracy, a good
choice of~$c$ (with the log correction term empirically determined;
see~\cite{lw} and~\cite[Thm.~2.4 and Prop.~17]{orx} for relevant theory) is
\begin{equation}
c = \pi n - 12 \log n \qquad (n\ge 8).
\label{c}
\end{equation}
Figure~\ref{ggqcos} shows the accuracy of the band-limited
quadrature formula with $n=50$ for approximating $\exp(ikx)$.
Figure~\ref{figtransggq} repeats the convergence comparison
of Figures~\ref{figcc1} and~\ref{figtrans}.\footnote{The
Chebfun calculations made use of
{\tt c = pi*n-12*log(n),}
{\tt [s,w] = pswfpts(n,c,\textquotesingle ggq\textquotesingle )}.}

The term ``band-limited quadrature'' is convenient, but
misleading.  
It is true that if one has an integrand that is
band-limited to a known range $k\in [-c,c\kern .5pt]$, or nearly so,
then these formulas
will be excellent.  But as equation (\ref{c}) and
Figure~\ref{figtransggq} illustrate,
this is not the only potential application, and to say that these formulas 
are for integration of band-limited functions is
not unlike saying that Gauss quadrature is for integration of polynomials.

Both of the alternatives to Gauss quadrature we 
have outlined in this section lead to the conclusion that the potential gain 
is a factor of $\pi/2$ in convergence rate.  This is hardly
enough to be very important in practice (at least in one space
dimension~\cite[Fig.~11]{cube}), and I share the view,
which is discussed with particular substance in section 6.9
of~\cite{bp}, that Gauss quadrature is the best choice known for
general use.

\section{\label{secgh}Gauss--Hermite quadrature}

Gauss--Hermite quadrature is the standard Gauss quadrature method
for integrals over the whole real line.
One supposes that a function $f$ is given and one wants to
compute the integral 
\begin{equation}
I = \int_{-\infty}^\infty e^{-x^2} \kern -2pt f(x) 
\kern .7pt dx,
\label{gaussherm}
\end{equation}
which is (\ref{I}) with $D = (-\infty,\infty)$ modified
by the introduction of the weight function $\exp(-x^2)$.
The approximation will take the form (\ref{In}) as usual,
now with nodes $\{x_j\}$ in $(-\infty,\infty)$,
and the Gauss--Hermite formula is defined by the
nodes and weights taking the unique values such
that $I(f) = I_n(f)$ whenever $f$ is a polynomial of
degree $2n-1$.  According to
Gautschi~\cite{gautschi}, this method was introduced by
Gourier in 1883~\cite{gourier}.
As with its unweighted progenitor discussed in the last section,
Gauss--Hermite quadrature is intimately associated with orthogonal
polynomials~\cite{szego}.  These are the {\em Hermite polynomials}
$H_0(x) = 1$, $H_1(x) = 2x$, $H_2(x) = 4x^2-2,\dots,$ which
are orthogonal over $(-\infty,\infty)$ with respect to the inner
product
\begin{equation}
\langle f,g\rangle = \int_{-\infty}^\infty
e^{-x^2} \kern 2pt \overline{f(x)} \kern 1pt g(x) \kern .7pt dx.
\end{equation}
The nodes $\{x_k\}$ are the roots of $H_n$, and as with
Gauss--Legendre quadrature, algorithms have been developed to
compute the nodes and weights with just $O(n)$ work~\cite{gst,tto}.
In Chebfun, {\tt [x,w] = hermpts(n)}.

In an application, one might start from an integrand $g(x)$
with Gaussian decay, so that $g$ can be written $g(x)
= \exp(-x^2) f(x)$ with $f(x) = O(1)$ on $(-\infty,\infty)$.
It has been recognized from the beginning that there may be
problems in practice with determining the right decay behavior.
What if $g(x)$ looks more like $\exp(-\sigma x^2)$ for some
$\sigma \ne 1$, or decays in a non-Gaussian manner?  We shall
stay away from these questions and assume that $\exp(-x^2)$
well approximates the shape of $g$.

Even in this most favorable setting, Gauss--Hermite quadrature is
terribly inefficient as $n\to\infty$.  For large $n$, most of the
nodes lie far enough out along the real axis that their weights
are minuscule\kern .5pt ; these terms then contribute negligibly to the
sum (\ref{In}) and can be thrown away!  This curious phenomenon
is known to experts, but I am unaware of any discussion in the
literature of where the reasoning that led to the Gauss--Hermite
formula has broken down.  How can a formula that is optimal in a
precise mathematical sense be so plainly suboptimal in practice?
We shall first illustrate the problem, and then give an answer
to this question.

\begin{figure}
\begin{center}
\vskip 8pt
\includegraphics[scale=.93]{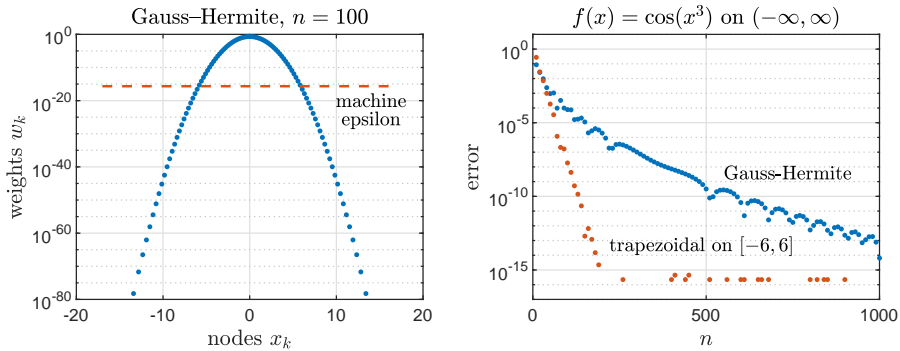}
\end{center}
\caption{\label{hermfig}A practitioner's view of the
inefficiency of Gauss--Hermite quadrature for large $n$.
Many of the weights $w_k$ lie below machine precision---a fraction
approaching $100\%$ as $n\to\infty$.  The convergence
rate is just root-exponential, whereas if one truncates to a finite
interval and applies a standard quadrature method, it becomes
exponential.  The missing dots for the trapezoidal rule
correspond to errors that are exactly zero in floating point
arithmetic.}
\end{figure}

Figure~\ref{hermfig} presents the inefficiency of Gauss--Hermite
quadrature as a practitioner might encounter it.  For the integrand
we take $f(x)=  \cos(x^3)$, mixing energy at all wave numbers.
The errors for Gauss--Hermite quadrature as a function of $n$
line up along a parabola on a semilog scale, corresponding to slow
root-exponential convergence, and even with $n=1000$ the error is
no smaller than $10^{-13}$.  Yet $\exp(-x^2)$ is of the order of
$10^{-16}$ or less for $|x|> 6$, so for practical purposes, this
integral might as well be posed on the compact interval $[-6,6]$.
Sure enough, the periodic trapezoidal rule applied on that interval
shows much faster exponential convergence down to machine precision, and the
convergence would be similarly fast for the Clenshaw--Curtis or
Gauss formulas on $[-6,6]$.  If $f(x)$ is changed to $\exp(-1/x^2)$
as in Figure~\ref{figcc1} (not shown), the convergence curves have
about the same shapes though with about half the convergence rate.

\begin{figure}
\begin{center}
\vskip 8pt
\includegraphics[scale=.81]{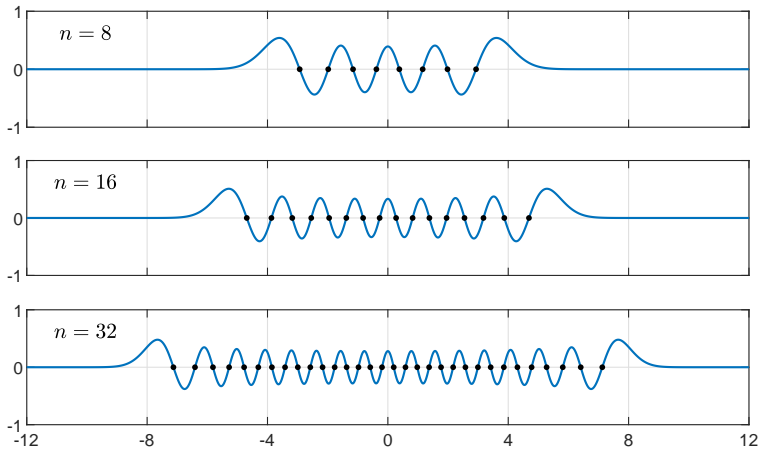}
\end{center}
\caption{\label{hermfig2}Hermite functions $(\ref{hermfuns})$,
a complete orthonormal set in $L^2(-\infty,\infty)$ whose
zeros are the Gauss--Hermite quadrature nodes $\{x_k\}$. The
interval around $x=0$ where the amplitude is $O(1)$ broadens
in proportion to $\sqrt n$, even though the integrand
of $(\ref{gaussherm})$ decays at the rate $\exp(-x^2)$.}
\end{figure}

The left side of Figure~\ref{hermfig} shows how the 
Gauss--Hermite formula wastes its effort.
For $n=100$, about half the weights
lie below machine precision, 48 of them to be exact,
and this fraction will increase with $n$.  For $n=1000$,
the number of weights below machine precision is 836.  Clearly
if $f$ is of order $1$, these function samples
will not contribute usefully to the evaluation of the integral.
Figure~\ref{hermfig2} illustrates the problem from another angle,
plotting the {\em Hermite functions}
\begin{equation}
\psi_n(x) = H_n(x)\kern .7pt e^{-x^2/2} (\sqrt {2\pi } \kern .7pt n!)^{-1/2}, 
\quad n = 0, 1, \dots
\label{hermfuns}
\end{equation}
for three values of $n$.  These functions form a complete orthonormal
set in the unweighted space $L^2(-\infty,\infty)$, and
Gauss--Hermite quadrature implicitly expands $f$ in this
basis.\footnote{Hermite functions have been well known to
physicists since the 1920s, for they are the eigenfunctions of the
Schr\"odinger equation for the harmonic oscillator.  In Chebfun,
try {\tt x = chebfun(\textquotesingle x\textquotesingle,[-3,3]),
quantumstates(x\^{~}\kern -4pt 2,25)}.} It is clear from the figure
that this is going to be be an inefficient way of capturing the
part of $f$ that matters to the integral, with $\psi_{32}$, for
example, taking values of size $O(1)$ all across the interval
$[-8,8\kern .8pt]$ even though $\exp(-8^2\kern .5pt) \approx
10^{-28}$.

Looking more deeply, there is no need for machine epsilon
to be part of the discussion.  For mathematical convergence
of $I_n(f)$ to $I(f)$ when $I_n$ is defined by quadrature
over a finite interval, the length of the interval will need
to grow as $n\to\infty$, and the right choice for analytic
functions $f(x)$ of size $O(1)$ on $(-\infty,\infty)$ is an
interval of size $[-O(n^{1/3}), O(n^{1/3})]$, which balances a
domain-truncation error of order $\exp(-n^{2/3})$ and a discretization
error of the same order since the sample step size will be $h
= O(n^{-2/3})$.  (Discretization errors for the trapezoidal
rule are quantified in~\cite{trap}; the tool for such estimates
is Cauchy integrals.)  From this balance one can quantify the
inefficiency of Gauss--Hermite quadrature, which places nodes
approximately uniformly in $[-O(n^{1/2}), O(n^{1/2})]$ whereas
$[-O(n^{1/3}), O(n^{1/3})]$ would suffice.  Thus in effect only
a fraction of order $n^{1/3}/n^{1/2} = n^{2/3}$ of the nodes are
utilized, meaning that Gauss--Hermite quadrature employs more nodes
than necessary by a factor of order $n^{1/3}$.  The ratio increases
to nearly order $n^{1/2}$ for nonanalytic functions $f$, where intervals
growing just logarithmically rather than algebraically with $n$ are
appropriate for balancing domain-truncation and discretization errors.

A theorem makes the point precise; the proof is given in the appendix.
\smallskip

\begin{theorem}
\label{thm}
Let $f(x)$ be analytic and bounded for $x\in (-\infty,\infty)$,
and suppose $\exp(-x^2)f(x)$ extends to a bounded analytic function
in an infinite strip $-a < \Im x < a$.
Let\/ $L>0$ be fixed, and for 
each $n\ge 1$, let\/ $I_n$ be the estimate of the integral\/
$I$ of $(\ref{gaussherm})$ obtained by applying 
Gauss--Legendre, Clenshaw--Curtis, or trapezoidal quadrature
on the truncated interval $[-Ln^{1/3}, L n^{1/3}\kern .7pt ]$.
Then for some
$C>0$,
\begin{equation}
|I-I_n| = O(\exp(-Cn^{2/3})), \qquad n\to \infty.
\label{est}
\end{equation}
\end{theorem}

In the literature, several authors have recommended adjustments
to Gauss--Her\-mite quadrature that are related in one manner
or another to truncating the domain to a finite interval such as
$[-O(n^{1/3}), O(n^{1/3})]$.  Mastroianni and Monegato and their
coauthors have published a number of papers in this direction,
focussing mainly on the analogous case of Gauss--Laguerre
quadrature on $[\kern .5pt 0,\infty)$~\cite{mm00,mm03}
(see next section).  Townsend,
Trogdon and Olver note that for large $n$, only about $25\sqrt
n$ of the weights $w_k$ are greater than the smallest machine
number ${\approx}\kern 1pt 10^{-308}$ in IEEE double precision,
and recommend ``subsampling'' to retain only nodes and weights
that will matter~\cite{tto}.  In unpublished work presented
at a conference in 2018, Weideman showed how
$O(\exp(-n^{2/3}))$ convergence with a particularly favorable
constant can be obtained by
applying the Gauss--Hermite formula with the weight function
$\exp(-(n/2)^{1/3} x^2)$ rather than $\exp(-x^2)$~\cite{weidtalk}.

\begin{figure}
\begin{center}
\vskip 10pt
\includegraphics[scale=.73]{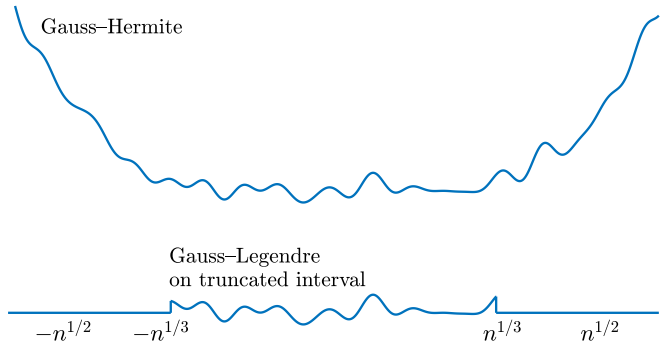}
\end{center}
\caption{\label{hermfig3}Schematic illustration of the approximation
problems $f\approx v\in V_n$ 
associated with Gauss--Hermite quadrature for $(\ref{gaussherm})$ and
the alternative method of applying the Gauss--Legendre formula
on a truncated interval $[-O(n^{1/3}),O(n^{1/3})]$, as in
Theorem~$\ref{thm}$.
The Gauss--Legendre method only needs to approximate $f$
on $[-O(n^{1/3}),O(n^{1/3})]$ to achieve an error $(\ref{1norm})$
of size $O(\exp(-O(n^{2/3})))$.  The Gauss--Hermite method is
effectively forced to approximate on a wider interval
$[-O(n^{1/2}),O(n^{1/2})]$ just to keep the growth of polynomials
under control.}
\end{figure}

However, the literature seems not to confront the conceptual
question:~what has gone wrong with the Gauss--Hermite notion
of optimality?  
Here is an answer,
summarized schematically in Figure~\ref{hermfig3}.
Suppose we have a quadrature formula for
(\ref{gaussherm}) that is exact
for all functions $v$ in an $n$-dimensional space $V_n$.
Its effectiveness will depend
on how well $f$ can be approximated by functions $v\in V_n$ in
the weighted $1$-norm
\begin{equation}
\|f-v\| = \int_{-\infty}^\infty e^{-x^2} |f(x) - v(x)|\kern .7pt dx.
\label{1norm}
\end{equation}
(The study of weighted approximation problems on the real line was
initiated by Bernstein~\cite{lub}.)
Now Gauss--Hermite quadrature corresponds to taking $V_n$
as the space of polynomials of degree $n-1$, and this space
is terribly inefficient for these weighted approximations
since polynomials grow so fast as $|x|$ increases.
For example, $x^n \exp(-x^2)$ reaches
a maximum of $(n/2\kern .3pt e)^{n/2}$ at $x= \sqrt{n/2}$ --- about
$10^6$ for $n=20$ and $10^{17}$ for $n=40$.
These huge numbers force a polynomial approximation to $f$
to pay attention to the whole interval
$[-O(n^{1/2}),O(n^{1/2})]$, just to keep (\ref{1norm}) under control.
By contrast, when we apply Gauss--Legendre quadrature on a
truncated interval $[-O(n^{1/3}),O(n^{1/3})]$, we are changing
the approximation space $V_n$ to the set of polynomials
of degree $n-1$ multiplied by the characteristic function of
this interval.  The approximation can focus
just on the short interval, where much better accuracy is achievable.

\section{\label{sec3}Three more examples}
We now mention briefly three further examples of quadrature
formulas for which exactness proves an inaccurate guide to accuracy.

\begin{figure}
\begin{center}
\vskip 8pt
\includegraphics[scale=.93]{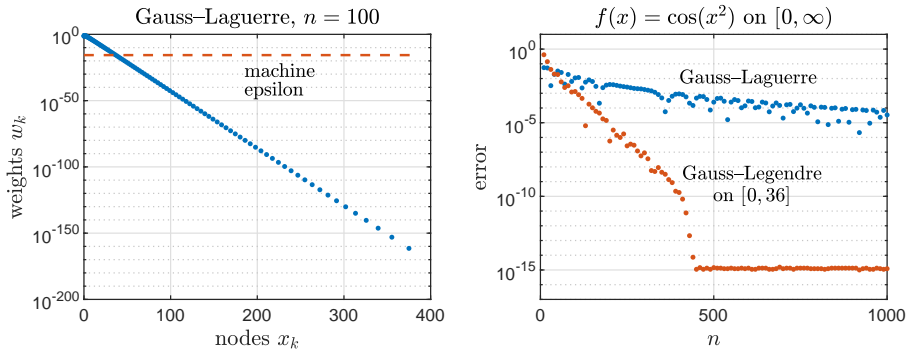}
\end{center}
\caption{\label{lagfig}Analogue of Figure~$\ref{hermfig}$ for
Gauss--Laguerre quadrature, now with $f(x) = \cos(x^2)$.  As
with Gauss--Hermite, many of the weights
lie below machine precision and the convergence
rate is just root-exponential; one does much better
by discarding nodes and weights or
truncating to a finite interval.}
\end{figure}

The first is Gauss--Laguerre quadrature, the analogue of 
(\ref{gaussherm}) for integration over a semi-infinite interval:
\begin{equation}
I = \int_0^\infty \kern -1pt e^{-x} f(x) 
\kern .7pt dx.
\label{gausslag}
\end{equation}
The design principle is that the nodes and weights should be
such that the formula is exact when $f$ is a polynomial of
degree $2n-1$.
(In Chebfun, they can be computed by {\tt [x,w] =
lagpts(n)}, an implementation of a fast algorithm of Huybrechs
and Opsomer~\cite{ho}.)
As with Gauss--Hermite quadrature, one
finds that many of the weights are so small that the corresponding
terms contribute negligibly to the result.
For $n=100$ the minimal
weight is $10^{-162}$ (decreasing exponentially with $n$), and
just 38 weights lie above standard machine precision (a fraction
increasing
as $O(n^{1/2})$).  For analytic $f$ of size $O(1)$ on $[\kern
.5pt 0,\infty)$ one has a domain-truncation error of size $\exp(-O(n))$
and a discretization error of order $\exp(-O(n^{1/2}))$, so the
two are out of balance.  By dropping nodes and weights or truncating
to a shorter interval one can do much better,
as illustrated in Figure~\ref{lagfig}.
Truncated Laguerre formulas have been investigated in detail
by Mastroianni and Monegato and their coauthors, though
mostly for more general functions $f$ with the property
that $e^{-x}f(x)$ may decay only algebraically as
$x\to\infty$~\cite{mm00,mm03}.

The second example involves the trapezoidal rule, that is,
the approximation of an integral of a function $f(x)$ by the
integral of its piecewise linear interpolant in a given set of
sample points.  Here the exactness principle is that $I_n(f) =
I(f)$ whenever $f$ is a piecewise linear function with breaks at
the nodes, which leads readily to an $O(h^2)$ accuracy bound for
any $f$ that is twice continuously differentiable.
However, there is an important special
case in which the convergence is much faster, which we exploited
in Figure~\ref{hermfig}: when $f$ is smooth
and periodic and the nodes are equally spaced.  If $f$ is analytic,
the convergence becomes exponential, at a rate $\exp(-O(n))$.
The exactness principle based on piecewise linears fails to detect
this, but it can be rescued (at least up to a Gauss factor of 2
attributable to aliasing) by the observation that for equispaced
points, the piecewise linear interpolant has the same integral as
a trigonometric interpolant.  (Proof: in both cases the integral
is equal to the mean of the sample data times the length of the
interval.)  For the trigonometric interpolant, fast convergence is
readily proved via contour integrals or aliasing of Fourier series.
Thus the special accuracy of the periodic equispaced trapezoidal
rule can be understood by an exactness principle after all,
so long as it is based on trigonometric interpolants instead of
piecewise linears.  For details see~\cite{trap}.

The third example concerns {\em cubature formulas\/} for
integration in the $s$-dimensional hypercube $D = [\kern .7pt
0,1]^s$.  Following an idea introduced by James Clerk Maxwell in
1877~\cite{maxwell}, one may approximate $I$ by the integral of
a degree $d$ multivariate polynomial interpolant through a set
of data values at $n$ points $x_k\in D$.  For each $d\ge 0$,
the dimension of the space of polynomials is $N = {s+d\kern
1pt\choose s}$, and one would hope to be able to have $n\le N$
nodes with positive weights.  A theorem of Tchakaloff asserts
that this is indeed possible (in greater generality, not just for
a uniformly weighted integrals over a hypercube), even though a
suitable set of nodes may be difficult to determine~\cite{tchak}.
For extensive information about cubature formulas, see the survey
article by Cools~\cite{cools}.

The difficulty with the exactness principle for cubature formulas
concerns the limitations of multivariate polynomials as a guide
to accuracy.  The standard definition of the degree $d$ of a
multivarate polynomial is the maximum 1-norm of its exponents; thus
$p(x,y) = x^3y^4$, for example, has degree $d=7$.  This definition
has the property of isotropy (i.e., rotation-invariance) in
$s$-space: if the variables are transformed by a rotation,
the degree of a polynomial $p$ is unchanged.  However, the
hypercube is itself far from isotropic.  The diagonals are
$\sqrt s$ times longer than the diameters, and there are $2^d$
of them, so most of its volume is ``in the corners'' in the
sense of lying outside the inscribed hypersphere.  As a consequence,
cubature formulas designed on Maxwell's principle will behave
anisotropically in $[\kern .7pt 0,1]^2$, giving better accuracy
for an integrand $f(x) = \varphi(x_1)$ aligned along one axis, say,
than for the same function rotated to $f(x) = \varphi((x_1+\cdots +
x_s)/\sqrt s\kern 1pt )$.  For angle-independent resolution in the
hypercube it would be necessary to base cubature formulas instead
on the {\em Euclidean degree}, defined in terms of the 2-norm of
the exponents.  (For example, the Euclidean degree of $x^3y^4$
is $5$.)  This effect was first pointed out in~\cite{cube}, and
a theorem making it precise was published in~\cite{hypercube}.
Note the analogy to the discussion of Gauss quadrature on $[-1,1]$
in section~\ref{secgauss}.
There the issue was translation-invariance, whereas
here it is rotation-invariance.

It would be interesting to provide a numerical illustration of
the suboptimality of standard cubature formulas.
In preparing this article, however, I have
come to realize that cubature formulas are not used much in
practice, making it unclear exactly what nodes and weights 
might be appropriate for
such a comparison.  It appears that most multiple integrals
are handled by essentially 1D methods such as tensor products, 
or by Monte Carlo methods and their relatives, or by other more
specialized tools~\cite{dks}.  Nevertheless we can estimate how much
efficiency should be lost, in principle, by building cubature formulas in the 
standard manner on the total degree.
The number of coefficients needed to specify a polynomial of
Euclidean degree $d$ is $d^s$ times the volume of an orthant of
the $s$-hypersphere,
$$
N_{\rm euclidean} = {d^s\pi^{s/2}\over 2^s (s/2)!}
\sim {d^s\over \sqrt{\pi s}\kern .7pt}
\left( {\pi e\over 2 s}\right)^{\kern -1.5pt s/2}\kern -2pt ,
$$
with the asymptotic approximation $\sim$ referring
to the limit $d\to\infty$.
To get the same resolution via the total degree
would require the number of coefficients to 
be $d^s$ times the volume of the $s$-simplex
expanded by $\sqrt s$ in each direction, 
$$
N_{\rm total} = {d^s s^{s/2}\over s!} \sim
{d^s e^s \over s^{s/2}\sqrt{2\pi s}\kern .7pt}
\kern 1pt.
$$
The inefficiency ratio associated with the standard design
principle of cubature formulas is accordingly
$$
{N_{\rm total} \over N_{\rm euclidean}} = {(s/2)!\over s!}
\left({4s\over\pi}\right)^{\kern -1.5pt s/2} \kern -2pt \sim
{1\over \sqrt 2 \kern .7pt } \left({2\kern .5pt e\over \pi}
\right)^{\kern -1.5pt s/2} \kern -2pt \approx (1.3)^{s-1}.
$$
For dimensions $s = 1$, $2$, $5$, $20$, and $40$ the ratios are about
$1$, $1.27$, $2.83$, $171$, and $41104$.

\section{\label{discussion}Discussion}
Quadrature theory is an edifice built over 200 years,
featuring both detailed estimates and general theories.  Often it
may be hard to extract the important points, not least because
valid estimates have a way of turning out to be far from sharp.
In this article I have focussed on the theme of the exactness
principle and how it may mislead, for this principle is the
basis of how quadrature is presented to students in our textbooks.
Unfortunately, on each particular topic, there are undoubtedly
relevant results in the literature I am unaware of,
for which I apologize.

Quadrature theory is not a hot research area nowadays; like complex
analysis, it is a field that we use all the time but which has a
way of seeming finished.  A book that I have found particularly
valuable is the 2011 monograph by Brass and Petras~\cite{bp}.
The opening chapter presents a ``standard estimation framework''
that elegantly makes precise the central question: given a
quadrature formula, how can we speak quantitatively of its degree
of optimality in integrating functions of particular classes?
The book is full of results and references to detailed work on
these problems.  A definite integral is just one example of a
linear functional, and the estimation framework for quadrature
is related to wide-ranging theories including $n$-widths,
information-based complexity, and optimal recovery.  A recent
contribution in this area is~\cite{dfpw}, and a textbook by Simon
Foucart is forthcoming that will present links to data science
and machine learning~\cite{foucart}.

The exactness principle for designing quadrature formulas is
algebraic, a matter of whether certain quantities are exactly
zero or not.  Quadrature, however, is a problem of analysis,
concerned with whether certain quantities are small or not.
It is to be expected that there will be some discrepancy
between the two.\footnote{See
the discussion of sampling theory and approximation theory on
p.~2118 of~\cite{austin}.}
Still, it is surprising that sometimes the discrepancy
can be huge without our having fully noticed it, as in
the case of Gauss--Hermite quadrature.

\section*{Appendix.  Proof of Theorem~\ref{thm}}

\indent~
\smallskip

\begin{proof}
Write $g(x) = \exp(-x^2) f(x)$, and for
fixed $L>0$ and each $n\ge 1$, let $J_n$ denote the
interval $J_n = [-Ln^{1/3},Ln^{1/3}\kern .7pt ]$.
Since $f(x)$ is bounded for $x\in (-\infty,\infty)$, truncating
$(-\infty,\infty)$ to $J_n$ introduces an
error in the integral of $g$
of size $O(\exp(-Cn^{2/3}))$ for some $C$.  Thus it is
enough to show that the Gauss--Legendre, Clenshaw--Curtis,
and trapezoidal approximations to the integral of $g$ over
$J_n$ also have accuracy of this order.

For the trapezoidal rule, this follows from contour integral
arguments as given in Section~5 of~\cite{trap}.  Theorem~5.1 of
that section has additional assumptions concerning integrals
of $g(x)$ along horizontal lines in the strip and decay of $g(x)$ to
$0$ as $|x|\to \infty$, but these are only needed to derive
an exact decay exponent.  For our purposes, where the claim is
just accuracy $O(\exp(-Cn^{2/3}))$ for some $C$, boundedness of $g$
is enough.

For the Gauss--Legendre and Clenshaw--Curtis formulas, the error
estimate follows from Theorem~19.3 of~\cite{atap}.
That theorem asserts $O(\kern .7pt \rho^{-n})$ accuracy for
integration of an analytic function on $[-1,1]$, assuming
it can be analytically continued to a bounded analytic function
in the Bernstein $\rho$-ellipse with foci $\pm 1$.
Here, the ellipse must be scaled to foci $\pm Ln^{1/3}$,
and for it to lie in the given fixed strip of analyticity as $n\to\infty$,
we can take decreasing values $\rho = 1 + d \kern .5pt n^{-1/3}$ for a
sufficiently small fixed value of $d$. 
The rescaled Theorem~19.3 of~\cite{atap} now gives errors of size
$O(n^{1/3}(1 + d\kern .5pt n^{-1/3})^{-n}) = O(\exp(-Cn^{2/3}))$, as required.
\end{proof}

\section*{Acknowledgments}
I am grateful for advice and suggestions from Jim Bremer, Ron
DeVore, Simon Foucart, Walter Gautschi, Andrew Gibbs,
Abi Gopal, Karlheinz Gr\"ochenig, Nick Hale, 
Daan Huybrechs, Giuseppe Mastroianni, Giovanni
Monegato, Yuji Nakatsukasa, Incoronata Notarangelo,
Sheehan Olver, Allan Pinkus, Kirill Serkh, Ian Sloan,
Alex Townsend, and Andr\'e Weideman.
Particularly important contributions were
made by Bremer and Hale in Section~4 (Gauss) and Weideman in
Section~5 (Gauss--Hermite).

\end{document}